\documentclass{birkjour}
\usepackage{amsmath}
\usepackage{amssymb}
\usepackage[all]{xypic}
\newtheorem{thm}{Theorem}[section]
\newtheorem{cor}[thm]{Corollary}
\newtheorem{lem}[thm]{Lemma}
\newtheorem{prop}[thm]{Proposition}
\theoremstyle{definition}
\newtheorem{defn}[thm]{Definition}
\theoremstyle{remark}
\newtheorem{rem}[thm]{Remark}
\newtheorem{ex}[thm]{Example}
\long\def\Thm#1{\begin{thm} #1 \end{thm}}
\long\def\Cor#1{\begin{cor} #1 \end{cor}}
\long\def\Lem#1{\begin{lem} #1 \end{lem}}
\long\def\Prop#1{\begin{prop} #1 \end{prop}}
\long\def\Def#1{\begin{defn} #1 \end{defn}}
\long\def\Rem#1{\begin{rem} #1 \end{rem}}
\long\def\Ex#1{\begin{ex} #1 \end{ex}}
\def\Sect{\section}
\def\imm{\looparrowright}
\def\Rarr#1#2{\xrightarrow[#2]{#1}}
\long\def\Ref#1#2#3#4#5#6{
\bibitem{#1}
{\rm #2,}
\textit{#3.}
{\rm #4}
\textbf{#5}
{\rm #6.}
}
\long\def\Refb#1#2#3#4{
\bibitem{#1}
{\rm #2,}
\textit{#3.}
#4.
}
\def\ZZ{{\Zz /2}}
\def\CC{{\mathcal C}}
\def\CCC{\mathfrak{L}(f)}
\def\Dp{\mathfrak{D}(f)}
\def\Dpone{\mathfrak{D}(f')}
\def\barDp{\overline{\mathfrak{D}}(f)}
\def\barDpone{\overline{\mathfrak{D}}(f')}
\def\hDp{\hbox{\rm h-${\mathfrak{D}(f)}$}}
\def\hI{\hbox{\rm h-${\mathfrak{I}(f^{(1)},f^{(2)})}$}}
\def\barhDp{\hbox{\rm h-$\overline{\mathfrak{D}}(f)$}}
\def\barD{Z}
\def\CD{\hbox{\rm h-$\widehat{\mathfrak{D}}(f)$}}
\def\CDa{\hbox{\rm h-$\widehat{\mathfrak{D}}(f^{(1)})$}}
\def\CDb{\hbox{\rm h-$\widehat{\mathfrak{D}}(f^{(2)})$}}
\def\CDone{\hbox{\rm h-$\widehat{\mathfrak{D}}(f')$}}

\def\EE{{\mathcal{D}}}
\def\O{{\rm O}}
\def\Zz{{\mathbb Z}}
\def\Rr{{\mathbb R}}
\def\Cc{{\mathbb C}}

\def\Hom{{\rm Hom}}

\def\into{\hookrightarrow}
\def\imm{\looparrowright} 
\def\iso{\cong}
\def\comp{\mathbin{\mathchoice
{\circ}
{{\scriptstyle\circ}}
{{\scriptscriptstyle\circ}}
{{\scriptscriptstyle\circ}}
}}
\def\st{\mid} 
\def\Null{{\rm Null}} 
\def\hNull{\hbox{\rm h-Null}} 
\def\hInd{\hbox{\rm h-$\gamma$}} 
\def\map{{\rm map}}
\def\geq{\geqslant} \def\leq{\leqslant} 
\begin{document}

\title{The geometric Hopf invariant and double points}
\author{Michael Crabb}
\address{%
Institute of Mathematics,\\
University of Aberdeen, \\
Aberdeen AB24 3UE\\
Scotland, UK}
\email{m.crabb@maths.abdn.ac.uk}
\author{Andrew Ranicki}
\address{%
School of Mathematics,\\
University of Edinburgh,\\
Edinburgh, EH9 3JZ\\
Scotland, UK
}
\email{a.ranicki@ed.ac.uk}

\begin{abstract}
The geometric Hopf invariant of a stable map $F$
is a stable $\ZZ$-equivariant map $h(F)$ such that the stable
$\ZZ$-equivariant homotopy class of $h(F)$ is the primary obstruction
to $F$ being homotopic to an unstable map. In this paper we express
the geometric Hopf invariant of the Umkehr map $F$ of an immersion
$f:M^m \imm N^n$ in terms of the double point set of $f$.
We interpret
the Smale-Hirsch-Haefliger regular homotopy classification of immersions
$f$ in the metastable dimension range $3m<2n-1$ (when a generic $f$ has no
triple points) in terms of the geometric Hopf invariant.
\end{abstract}

\maketitle
17 May 2010
%

\Sect*{Introduction}

The original Hopf invariant $H(F) \in \Zz$ of a map $F:S^3 \to S^2$
was interpreted by
Steenrod as the evaluation of the cup product in the mapping cone $C(F)$.
The mod 2 Hopf invariant $H_2(F) \in \ZZ$ of a map $F:S^j \to S^k$ was then
defined using the functional Steenrod squares of $F$.
The geometric Hopf invariant of a stable map $F:\Sigma^{\infty}X \to
\Sigma^{\infty}Y$ is the stable $\ZZ$-equivariant map
$$
h(F)~=~(F\wedge F)\Delta_X-\Delta_YF~:~\Sigma^\infty X\to
\Sigma^\infty(Y\wedge Y)
$$
measuring the failure of $F$ to preserve the diagonal maps of $X$ and
$Y$, with $\ZZ$ acting by the identity on $X$ and by the
transposition $T:(y_1,y_2) \mapsto (y_2,y_1)$
on $Y \wedge Y$. Thus $h(F)$ is a homotopy-theoretic generalization of the
functional Steenrod squares. The stable homotopy
class of $h(F)$ is the primary obstruction to $F$ being homotopic to an
unstable map.

Given an immersion $f:M^m \imm N^n$ with normal bundle $\nu (f)$ we
express the geometric Hopf invariant $h(F)$ of the Umkehr map
$F:\Sigma^{\infty}N^+\to \Sigma^{\infty}M^{\nu(f)}$ in terms of the
double point set of $f$, where $M^{\nu(f)}$ denotes the Thom space.
There are many antecedents for this
expression! The stable homotopy class of $h(F)$ depends only on the
regular homotopy class of $f$.  If $f$ is regular homotopic to an
embedding then $h(F)$ is stably null-homotopic.  We interpret the
Smale-Hirsch-Haefliger regular homotopy classification of immersions
$f$ in the metastable dimension range $3m<2n-1$ (when a generic $f$ has
no triple points) in terms of the geometric Hopf invariant.

In \cite{Hopf} we shall provide a considerably more detailed exposition
of the geometric Hopf invariant $h(F)$ and its applications to
manifolds.  This will include the $\pi_1(N)$-equivariant geometric Hopf
invariant $\widetilde{h}(F)$ needed for a homotopy-theoretic treatment
of the double point invariant $\mu(f)$ of Wall \cite{Wall} for a
generic immersion $f:M^m \imm N^{2m}$ which plays such an important
r\^ole in non-simply-connected surgery theory, with $M=S^m$.  When both
$M$ and $N$ are connected and oriented and $f$ induces the trivial map
$\pi_1(M)\to \pi_1(N)$, $\mu (f)$ is an element of the group
$$
\Zz[\pi_1(N)]/\langle g-(-)^mg^{-1}\,\vert\, g \in \pi_1(N)\rangle
$$
and $\mu(f)=0$ if (and, for $m > 2$, only if)  $f$ is regular
homotopic to an embedding, by the Whitney trick for removing double points.
In \cite{Hopf} $\widetilde{h}(F)$ will be shown to induce the
quadratic construction $\psi_F$ of \cite{Ranicki} on the chain level.

The present paper is set out as follows.  Section 1 describes briefly
the construction of the geometric Hopf invariant and its fibrewise
generalization.  The double point theorem is stated and proved in
Section 2.  In Section 3, building on work of Dax \cite{Dax}, Hatcher
and Quinn \cite{Quinn}, Salomonsen \cite{Salomonsen} and Li, Liu and
Zhang \cite{Ping}, we relate the geometric Hopf invariant in a stable
range to Haefliger's obstruction to the existence of a regular homotopy
from an immersion to an embedding.  The papers of Boardman and Steer
\cite{BoardmanSteer} and Koschorke and Sanderson \cite{KoschorkeSanderson}
are also relevant. The variation of the geometric Hopf
invariant of an immersion under a (not necessarily regular) homotopy is
computed in Section 4 in terms of the Smale-Hirsch-Haefliger
classification.  In Section 5 we use Whitney's figure-of-eight
immersion \cite{Whitney} to construct, in a metastable range,
immersions close to a given embedding.  Prerequisites for that section,
on the differential-topological classification of vector bundle
monomorphisms, are given in an Appendix.

We shall write the one-point compactification of a locally compact
Hausdorff topological space $X$ as $X^+$.
A subscript `$+$' will be used for the adjunction of a disjoint
basepoint to a space. If $X$ is compact $X^+=X_+$.
For a Euclidean vector bundle $\xi$ over a general space $X$,
we write $D(\xi )$ for the closed
unit disc bundle, $S(\xi )$ for the sphere bundle and $B(\xi )$ for
open unit disc bundle. The Thom space of $\xi$ is written as
$X^\xi$.  To simplify notation, we shall sometimes write
$Y^\xi$, rather than $Y^{p^*\xi}$, for the Thom space of the
pullback $p^*\xi$ by a map $p:Y\to X$, if the map $p$ is clear from
the context. Similarly, we sometimes write $V$, instead of $X\times V$,
for the trivial vector bundle over $X$ with fibre the vector space $V$.

Methods of fibrewise homotopy theory will be used extensively.
Fibrewise constructions, such as the one-point compactification,
the Thom space, or the smash product, over a base $B$ will be indicated
by attaching a subscript `$B$' to the relevant symbol.
We follow  the notation for (fibrewise) stable homotopy adopted in
\cite{coin}.
Consider fibrewise pointed spaces $X\to B$ and $Y\to B$ over an
ENR base $B$.
If $B$ is compact and $A$ is a closed sub-ENR, we write
$$
\omega^0_B\{ X;\, Y\} \text{\quad and\quad}
\omega^0_{(B,A)}\{ X;\, Y\},
$$
respectively,
for the abelian group of stable fibrewise maps $X\to Y$ over $B$
and the relative group
defined in terms of homotopy classes of maps that are zero over the
subspace $A$.
(See, for example, \cite[Part II, Section 3]{FHT}.)
We also need to consider fibrewise maps with compact supports.
When $B$ is not necessarily compact we write
$$
{}_c\omega^0_B\{ X;\, Y\}
$$
for the group of fibrewise stable maps that are zero outside a
compact subspace of $B$. The $\omega^0$-theories are extended,
using the fibrewise suspension $\Sigma_B$ over $B$, to
$\omega^i$-cohomology theories indexed by $i\in\Zz$. When $Y\to B$
is a trivial bundle $B \times S^i \to B$,
there are natural identifications of the
fibrewise groups with the reduced stable cohomotopy of an
appropriate pointed space:
$$
\omega^0_B\{ X;\, B\times S^i\} =\tilde\omega^i(X/B),
\text{\quad and\quad}
\omega^0_{(B,A)}\{ X;\, B\times S^i\} =\tilde\omega^i(X/(X_A\cup B)),
$$
where $X_A\to A$ denotes the restriction of $X\to B$.

We shall also need $\ZZ$-equivariant stable homotopy theory, which we
indicate by a prefix as ${}^\ZZ\omega^*$. So, for example, the
equivariant stable cohomotopy of a point, ${}^\ZZ\omega^0(*)$, is
the direct limit over all finite-dimensional $\Rr$-vector spaces $V$ and $W$
of the homotopy classes of pointed $\ZZ$-maps
$(V\oplus LW)^+\to (V\oplus LW)^+$.
Here, and throughout the paper, we write $L$ for  the non-trivial
$1$-dimensional representation $\Rr$ of $\ZZ$ with
the involution $-1$, and,
for a finite-dimensional real vector space $W$, abbreviate the tensor product
$L\otimes W$ to $LW$.

\smallskip

We thank Mark Grant for pointing out the relevance of \cite{Quinn}.
\Sect{A review of the geometric Hopf invariant}
Let $X$ and $Y$ be pointed topological spaces. It is convenient to assume that $X$ is
a compact ENR and that $Y$ is an ANR.

We shall be working both with the unreduced suspension of an unpointed
space $A$
$$s A ~=~[0,1] \times A/\!\sim~,~(0,a) \sim (0,a')~,~(1,a) \sim (1,a')~
(a,a' \in A)~,$$
and with the reduced suspension of a pointed space $(B,* \in B)$
$$\Sigma B~=~ [0,1] \times B/\!\sim~,~(0,b) \sim (1,b')~\sim (t,*)~
(b,b' \in B,\,t \in [0,1])~.$$
\indent Let $V$ be a finite dimensional Euclidean space, and let
$$S(V)~=~\{u \in V\,\vert\, \Vert u \Vert = 1\}$$
be the unit sphere in $V$ with respect to an inner product
$\Vert~\Vert$ on $V$. For $t \in [0,1]$, $u \in S(V)$ we write
$$[t,u]~=~\dfrac{tu}{1-t} \in V^+$$
with $[0,u]=0$, $[1,u]=+ \in V^+$. The maps
$$\begin{array}{l}
s S(V) \to V^+~;~(t,u) \mapsto [t,u]~,\\[1ex]
\Sigma S(V)_+ \to V^+/\{0,+\}~;~(t,u) \mapsto [t,u]
\end{array}
$$
are homeomorphisms.

Let the generator $T \in \ZZ$ act on $X \wedge X$ by transposition
$$T~:~X \wedge X \to X \wedge X~;~(x,y)\mapsto (y,x)\, .$$
The diagonal map
$$\Delta_X~:~X\to X \wedge X~;~x \mapsto (x,x)$$
is $\ZZ$-equivariant. The diagonal map $\Delta_{V^+ }$
extends to a $\ZZ$-equivariant homeomorphism
$$\kappa_V~:~LV^+ \wedge V^+ \to V^+ \wedge V^+~;~
(u,v) \mapsto (u+v,-u+v)\, .$$
The $\ZZ$-action $u \mapsto -u$ on $LV$ has fixed point $\{0\}$; the
$\ZZ$-action on the unit sphere $S(LV)$
is free.

The geometric Hopf invariant of a map
$F:V^+ \wedge X \to V^+ \wedge Y$
measures the difference $(F\wedge F)\Delta_X - \Delta_YF$, given that
$(F\wedge F)\Delta_{V^+ \wedge X}=\Delta_{V^+ \wedge Y}F$.
The diagram of $\ZZ$-equivariant maps
$$
\xymatrix@R+20pt@C+40pt{
LV^+ \wedge V^+ \wedge X \ar[r]^-{\displaystyle{1 \wedge \Delta_X}}
\ar[d]_-{\displaystyle{1 \wedge F}} &
LV^+ \wedge V^+ \wedge X \wedge X
\ar[d]^-{\displaystyle{G}}\\
LV^+ \wedge V^+ \wedge Y
\ar[r]^-{\displaystyle{1 \wedge \Delta_Y}}
&LV^+ \wedge V^+ \wedge Y \wedge Y}
$$
does not commute in general, with $G$ defined by
$$
\begin{array}{l}
G~=~(\kappa^{-1}_V \wedge 1)(F \wedge F)(\kappa_V \wedge 1)~:\\[1ex]
\hskip50pt
LV^+ \wedge V^+ \wedge X \wedge X
\to LV^+ \wedge V^+ \wedge Y \wedge Y~;\\[1ex]
\hskip75pt (u,v,x_1,x_2) \mapsto ((w_1-w_2)/2,(w_1+w_2)/2,y_1,y_2)\\[1ex]
\hskip75pt (F(u+v,x_1)=(w_1,y_1),F(-u+v,x_2)=(w_2,y_2))\, .
\end{array}
$$
However, the $\ZZ$-equivariant maps defined by
$$
\begin{array}{l}
p~=~G(1\wedge \Delta_X)~:~LV^+ \wedge V^+ \wedge X \to
LV^+ \wedge V^+ \wedge Y \wedge Y~;\\[1ex]
\hskip75pt (u,v,x) \mapsto ((w_1-w_2)/2,(w_1+w_2)/2,y_1,y_2)\\[1ex]
\hskip75pt (F(u+v,x)=(w_1,y_1),F(-u+v,x)=(w_2,y_2))\, ,\\[1ex]
q~=~(1 \wedge \Delta_Y)(1 \wedge F)~:~
LV^+ \wedge V^+ \wedge X \to
LV^+ \wedge V^+ \wedge Y \wedge Y~;\\[1ex]
\hskip75pt (u,v,x) \mapsto (u,w,y,y)~~(F(v,x)=(w,y))
\end{array}
$$
agree on $0^+ \wedge V^+ \wedge X =V^+ \wedge X \subseteq
LV^+ \wedge V^+ \wedge X$, with
$$
\begin{array}{l}
p\vert~=~q\vert~=~(\kappa^{-1}_V \wedge 1)\Delta_{V^+ \wedge Y}
F~=~(\kappa^{-1}_V \wedge 1)(F \wedge F)(\kappa_V \wedge 1)
\Delta_{V^+ \wedge X}~:\\[1ex]
\hskip100pt V^+ \wedge X
\to LV^+ \wedge V^+ \wedge Y\wedge Y\, .
\end{array}
$$

\Def{(\cite[pp.~306--308]{FHT})
The {\it geometric Hopf invariant} of a pointed
map $F:V^+ \wedge X \to
V^+ \wedge Y$ is the $\ZZ$-equivariant map given by the relative
difference of the $\ZZ$-equivariant maps $p,q$
\index{geometric Hopf invariant, $h_V(F)$}
$$
\begin{array}{ll}
h_V(F)~=~\delta(p,q)~:&\Sigma S(LV)_+ \wedge V^+ \wedge X \to
LV^+ \wedge V^+ \wedge Y\wedge Y~;\\[1ex]
&(t,u,v,x) \mapsto \begin{cases}
q([ 1-2t,u],v,x)&\hbox{if $0 \leq  t \leq  1/2$}\\[1ex]
p([ 2t-1,u],v,x)&\hbox{if $1/2 \leq  t
\leq 1$}
\end{cases}\\[3ex]
&(t \in [0,1],\, u \in S(LV),\, v \in V,\, x \in X)\, .
\end{array}
$$
}

We are primarily interested in the stable $\ZZ$-equivariant class
of $h_V(F)$, which depends only on the homotopy class of $F$,
and which we write simply as
$$
h_V(F) \in {}^\ZZ\omega^0\{ \Sigma S(LV)_+\wedge X;\,
LV^+\wedge (Y\wedge Y)\}\, .
$$
Using duality, we can rewrite this stable homotopy group in different ways and
thus obtain two other versions of the geometric Hopf invariant as follows.
Smashing $h_V(F)$ with the identity on the sphere $S(LV)_+$
and composing with the duality map
$$
LV^+ \to \Sigma S(LV)_+ \to
S(LV)_+\wedge \Sigma S(LV)_+
$$
we get a map
$$
V^+ \wedge LV^+ \wedge X \to
V^+ \wedge LV^+ \wedge S(LV)_+ \wedge (Y\wedge Y)
$$
and a second version of the geometric Hopf invariant as an element
$$
h'_V(F)\in{}^\ZZ\omega^0\{ X;\, S(LV)_+\wedge (Y\wedge Y)\}\, .
$$
\Rem{\label{remark1}
The $\ZZ$-equivariant cofibration sequence
$$
S(LV)_+ \to S^0=\{0\}^+ \to LV^+
$$
induces an exact sequence of stable $\ZZ$-equivariant homotopy groups
$$
{}^\ZZ\omega^0\{ X;\, S(LV)_+\wedge (Y\wedge Y)\} \to
{}^\ZZ\omega^0\{ X;\, Y\wedge Y\} \to
{}^\ZZ\omega^0\{X ;\, LV^+\wedge Y\wedge Y\}\, .
$$
The stable class $h'_V(F)$ has image
$(F \wedge F)\Delta_X - \Delta_YF$ in
${}^\ZZ\omega^0\{ X;\, Y\wedge Y\}$.
}

We may also use the Adams isomorphism \cite[Theorem 5.3]{adams}
$$
{}^\ZZ\omega^0\{ X;\, S(LV)_+\wedge (Y\wedge Y)\}~\cong~
\omega^0\{ X;\, S(LV)_+\wedge_{\ZZ} (Y\wedge Y)\}
$$
to regard $h'_V(F)$ as a non-equivariant class
$$
h''_V(F)\in\omega^0\{ X;\, S(LV)_+\wedge_{\ZZ} (Y\wedge Y)\}\, .
$$
\Rem{Instead of applying the Adams isomorphism,
we can use fibrewise techniques.
The unstable map $h_V(F)$ lifts to
an equivariant fibrewise pointed map over $S(LV)$:
$$
S(LV)\times (\Sigma V^+ \wedge X) \to
S(LV)\times (V^+\wedge LV^+ \wedge  (Y\wedge Y))\, ,
$$
and we may divide by the free $\ZZ$-action to get
a fibrewise pointed map over the real projective space $P(V)=S(LV)/\ZZ$:
$$
P(V)\times (\Sigma V^+ \wedge X) \to
(V\oplus (H\otimes V))^+_{P(V)}\wedge_{P(V)}
(S(LV)\times_{\ZZ} (Y\wedge Y))\, ,
$$
where $H=S(LV)\times_{\ZZ} L$ is the Hopf line bundle over $P(V)$
and $\{-\}^+_{P(V)}$ denotes fibrewise one-point
compactification over $P(V)$.
This fibrewise map represents a stable class in the group
$$
\omega^{-1}_{P(V)}\{ P(V)\times X;\, (H\otimes V)^+_{P(V)}
\wedge_{P(V)} (S(LV)\times_{\ZZ}(Y\wedge Y))\}\, ,
$$
which may be identified by fibrewise Poincar\'e-Atiyah duality
(since $\Rr\oplus\tau P(V)= H\otimes V$) with
$$
\omega^0\{ X;\, S(LV)_+\wedge_{\ZZ} (Y\wedge Y)\}\, .
$$
(For the duality theorem, see, for example, \cite[Part II, Section 12]{FHT}.)
}
\Rem{By \cite[pp.~60--61]{Z2} the limit over all finite-dimensional
$V$ of the exact sequences in Remark \ref{remark1} is
a split exact sequence
$$
\omega^0\{ X;\, (E\ZZ )_+\wedge_{\ZZ} (Y\wedge Y)\} \Rarr{\delta}{}
{}^\ZZ\omega^0\{ X;\, Y\wedge Y\} \Rarr{\rho}{}\omega^0\{X ;Y\}
$$
with the fixed point map $\rho$ split by
$$
\omega^0\{X ;Y\} \to {}^\ZZ\omega^0\{ X;\,
Y\wedge Y\}~;~F \mapsto \Delta_YF\, .
$$
The stable $\ZZ$-equivariant homotopy class
$$
h(F)~=~(F\wedge F)\Delta_X-\Delta_YF \in {}^\ZZ\omega^0\{ X;\, Y\wedge Y\}
$$
is the image under $\delta$ of
$$
h'(F)\,=\,\varinjlim\limits_V \,  h'_V(F) \in
\omega^0\{ X;\, (E\ZZ )_+\wedge_{\ZZ} (Y\wedge Y)\}\, .
$$
}

We shall need some fundamental properties of the geometric Hopf invariant.
First, it is an obstruction to desuspension.
The Hopf invariant of a suspension is zero, but there is a more
precise result for the geometric Hopf invariant.

Given a pointed space $(A,* \in A)$ and $a\neq * \in A$ let
$a_1,a_2 \in A \vee A$ denote the two distinct copies of $a$,
and write
$$A \vee_a A~=~(A \vee A)/(a_1 \sim a_2)~.$$
\Lem{\label{suspension}
{\rm (i)} The $\ZZ$-equivariant maps
$$\begin{array}{l}
i_V~:~S(LV) \to \Sigma S(LV)_+~;~v \mapsto (1/2,v)~,\\[1ex]
j_V~:~\Sigma S(LV)_+\to LV^+\vee_0 LV^+~;~
(t,u) \mapsto \begin{cases}
[1-2t,u]_1&\hbox{if $0 \leqslant t \leqslant 1/2$}\\
[2t-1,u]_2&\hbox{if $1/2 \leqslant t \leqslant 1$}~,
\end{cases}\\[1ex]
k_V~:~LV^+\vee_0 LV^+ \to LV^+~;~v_1 \mapsto v,~v_2 \mapsto v
\end{array}$$
fit into a $\ZZ$-equivariant cofibration
$$\xymatrix{S(LV) \ar[r]^-{i_V} &\Sigma S(LV)_+ \ar[r]^-{j_V} & LV^+\vee_0 LV^+
\ar[r]^-{k_V}& LV^+}$$
with the composite $k_V \circ j_V:\Sigma S(LV)_+\to LV^+$
a $\ZZ$-equivariantly null-homotopic map.\\
{\rm (ii)} The geometric Hopf invariant of
$F:V^+ \wedge X \to V^+\wedge Y$ is the composite
$$\begin{array}{l}
h_V(F)~=~(q \vee_0 p)(j_V \wedge 1_{V^+ \wedge X})~:\\[1ex]
\Sigma S(LV)_+ \wedge V^+ \wedge X \to
(LV^+\vee_0LV^+)\wedge V^+ \wedge X \to LV^+\wedge V^+ \wedge Y \wedge Y
\end{array}$$
{\rm (iii)} Suppose that $F$ is the suspension $1_{V^+}\wedge F_0$
of a map $F_0 : X\to Y$, so that
$$p~=~q~=~1 \wedge (\Delta_Y \circ F_0)~:~LV^+ \wedge V^+ \wedge X \to
LV^+ \wedge V^+ \wedge Y \wedge Y~.$$
Then
$$h_V(F)~=~1 \wedge (k_V \circ j_V) \wedge (\Delta_Y\circ F_0)~ \simeq~ *$$
is $\ZZ$-equivariantly null-homotopic.
}
The second property, giving a formula for the Hopf invariant of a composition,
is suggested by the stable form described in Remark \ref{remark1}.
\Prop{\label{composition}{\rm (Composition formula).}
Let $X,\, Y$ and $Z$ be pointed spaces, and
suppose that $F:V^+\wedge X\to V^+\wedge Y$ and $G: V^+\wedge Y\to V^+\wedge Z$
are pointed maps. Then
$$
h_V(G\comp F) = h_V(G)[F] + [G\wedge G]h_V(F)\in
{}^\ZZ\omega^0\{ \Sigma S(LV)_+\wedge X;\, LV^+\wedge (Z\wedge Z)\},
$$
where
$$
[F]\in {}^\ZZ \omega^0\{ X;\, Y\} \text{\quad and \quad}
[G\wedge G]\in {}^\ZZ\omega^0\{ Y\wedge Y;\, Z\wedge Z\}
$$
are the stable classes determined by $F$ (with the trivial action of
$\ZZ$) and $G\wedge G$.
}
Lastly, the Hopf invariant satisfies the following sum formula.
\Prop{\label{sum}{\rm (Sum formula).}
Let $F_+,\, F_-$ be maps $V^+\wedge X\to V^+\wedge Y$.
Suppose that $v\in S(V)$. Choose a tubular neighbourhood of
$\{ v,\, -v\}$ in $V$ and let $\nabla : V^+\to V^+\wedge V^+$ be the
associated Pontryagin-Thom map.
Then the Hopf invariant of the sum $F=(F_+\wedge F_-)\circ (\nabla\wedge 1_X)$
is given by
$$
h_V(F)=h_V(F_+)+h_V(F_-) + \iota [(F_+\wedge F_-)\circ\Delta_X]
$$
where the induction homomorphism $\iota$ is the composition of the isomorphism
$$
\omega^0\{ X;\, Y\wedge Y\}
\iso
{}^\ZZ\omega^0\{ \Sigma S(Lv)_+\wedge X; \, (Lv)^+\wedge (Y\wedge Y)\}
$$
and the map
$$
{}^\ZZ\omega^0\{ \Sigma S(Lv)_+\wedge X;
\, (Lv)^+\wedge (Y\wedge Y)\}
\to {}^\ZZ\omega^0\{ \Sigma S(LV)_+\wedge X;\, LV^+\wedge (Y\wedge Y)\}
$$
induced by the inclusion of the $1$-dimensional subspace $\Rr v\into V$.
}
The explicit construction of the geometric Hopf invariant is readily
extended to the fibrewise theory.
Suppose now that $X\to B$ and $Y\to B$ are fibrewise pointed spaces
over an ENR $B$.  (We shall assume that both are locally fibre homotopy
trivial and that the fibres have the homotopy type of CW complexes,
finite complexes in the case of $X$.)
Consider a fibrewise pointed map $F: (B\times V^+)\wedge_B X \to
(B\times V^+)\wedge_B Y$. If $B$ is compact, we have a fibrewise
geometric Hopf invariant
$$
h_V(F) \in  {}^\ZZ\omega^0_B\{ (B\times \Sigma S(LV)_+)\wedge_B X;\,
(B\times LV^+)\wedge_B (Y\wedge_B Y)\}\, ,
$$
and corresponding variants $h'_V(F)$ and $h''_V(F)$.
See, for example, \cite[Part II, Section 14]{FHT}.
This fibrewise Hopf invariant is an obstruction to fibrewise desuspension.
Indeed, suppose that the restriction of $F$ to a closed sub-ENR $A\subseteq
B$ is the fibrewise suspension of a map $X_A \to Y_A$ over $A$.
Then Lemma \ref{suspension} (iii) shows how to define a relative fibrewise
Hopf invariant in
$$
{}^\ZZ\omega^0_{(B,A)}\{ (B\times \Sigma S(LV)_+)\wedge_B X;\,
(B\times LV^+)\wedge_B (Y\wedge_B Y)\}\, ,
$$
which lifts $h_V(F)$. (One uses the fact that the inclusion of $A$ in $B$
is a cofibration.)
When $B$ is not (necessarily) compact and $F$ is a fibrewise suspension
outside some compact subspace of $B$, the same method gives
a fibrewise Hopf invariant with compact supports:
$$
h_V(F) \in
{}_{\phantom{Z/}c}^\ZZ\omega^0_B\{ (B\times \Sigma S(LV)_+)\wedge_B X;\,
(B\times LV^+)\wedge_B (Y\wedge_B Y)\}\, .
$$

\Sect{The double point theorem} Let $f: M\imm N$ be a (smooth)
immersion of a closed manifold $M$ in a connected manifold $N$
(without boundary) of dimension $n$, with normal bundle $\nu (f)$,
usually abbreviated to $\nu$. We do not require $M$ to be
connected, nor that all the components should have the same
dimension; the maximum dimension of a component is denoted by $m$.
Let $e: M\to V$ be a smooth map to a finite-dimensional Euclidean
space $V$ such that $e(x)\not=e(y)$ whenever $f(x)=f(y)$ for
$x\not=y$. This gives a (smooth) embedding $(e,f): M \into V\times
N$ with normal bundle $V\oplus\nu$. (To be precise, there is a
short exact sequence $0\to V\to \nu (e,f)\to \nu \to 0$ which is
split by a choice of metrics.)

We have an associated Pontryagin-Thom map (defined up to homotopy)
$$
F : V^+\wedge N^+ \to V^+\wedge M^\nu .
$$
Its geometric Hopf invariant is a stable $\ZZ$-homotopy class
$$
h_V(F)\in{}^\ZZ\omega^0\{ \Sigma S(LV)_+ \wedge N^+  ;\,
 LV^+ \wedge (M^\nu\wedge M^\nu )\}\, ,
$$
where $\ZZ$ interchanges the factors of $M^\nu\wedge M^\nu$.

Suppose that the immersion $f$ is self-transverse and that
there are no $k$-tuple points for $k>2$.
(This is the case for a generic immersion if $3m<2n$.)
The {\it double point set}
$$
\Dp = \{ (x,y)\in M\times M-\Delta (M) \st f(x)=f(y)\}
$$
is then a smooth $\ZZ$-submanifold of $M\times M$
(of constant dimension $2m-n$ if $M$ is connected), on which
$\ZZ$ acts freely, and its normal bundle may be identified
with the pullback $j^*\tau N$ of
the tangent bundle of $N$ by the map $j:\Dp \to N$ mapping $(x,y)$ to
$f(x)=f(y)\in N$.
We also have a $\ZZ$-map $d :\Dp \to LV-\{ 0\}$ given
by $d(x,y)=e(x)-e(y)$, and thus an embedding
$$
(d,j) :  \Dp \into (LV-\{ 0\})\times N\, ,
$$
with normal bundle $LV\oplus k^*(\nu \times \nu )$,
where $k: \Dp\to M\times M$ is the inclusion.
The Pontryagin-Thom construction gives a $\ZZ$-map
$$
(LV-\{ 0\})^+\wedge N^+ \to
LV^+\wedge \Dp^{k^*(\nu\times\nu )}.
$$
Composing with the map induced by $k$, we get a $\ZZ$-homotopy
class
$$
\phi :\Sigma S(LV)_+ \wedge N^+ \to LV^+\wedge
(M^\nu\wedge M^\nu )\, .
$$
\Thm{\label{dpt}
{\rm (The double point theorem).}
The geometric Hopf invariant $h_V(F)$ of the Pontryagin-Thom map
$F$ is equal to the $\ZZ$-equivariant
stable homotopy class of the map $\phi$ determined, as described above,
by the double point manifold $\Dp$.
}
We may also consider the second version of the Hopf invariant
$$
h'_V(F) \in{}^\ZZ\omega^0\{ N^+ ;\,
S(LV)_+ \wedge (M^\nu\wedge M^\nu )\}\, .
$$
This, too, may be described directly in terms of the double points.
The embedding $\Dp \into LV\times N$ with normal
bundle $LV\oplus k^*(\nu\times\nu )$
provides a Pontryagin-Thom map
$$
LV^+\wedge N^+ \to
LV^+\wedge \Dp^{k^*(\nu\times\nu )}
$$
which we compose with the map $\Dp \to S(LV)\times (M\times M)$
given by $e$ and the inclusion $k$ to get
$$
\phi' : LV^+\wedge N^+ \to LV^+\wedge S(LV)_+
\wedge (M^\nu\wedge M^\nu )\, .
$$
\Cor{The second version $h'_V(F)$ of the geometric Hopf invariant
is represented by the $\ZZ$-map $\phi'$ defined in the text.
}
We also have the non-equivariant stable Hopf invariant
$$
h''_V(F) \in\omega^0\{ N^+ ;\,
(S(LV) \times_{\ZZ} (M\times M))^{\nu\times\nu}\}\, .
$$
The free $\ZZ$-manifold $\Dp$ is a double cover of the set
$\barDp\subseteq N$ of double points of the immersion $f$.
We have an induced map
$$
\barDp =\Dp /\ZZ \to S(LV)\times_{\ZZ} (M\times M)
$$
and an embedding of $\barDp$ in $N$ with normal bundle the pullback of
$\nu\times\nu$. The Pontryagin-Thom construction gives a map
$$
\phi'' :  N^+ \to \barDp^{\nu\times\nu} \to
(S(LV) \times_{\ZZ} (M\times M))^{\nu\times\nu}.
$$
\Cor{The stable geometric Hopf invariant $h''_V(F)$ is equal to the
stabilization $[\phi'']$ of the class determined by the double point
manifold $\barDp\subseteq N$.
}

These results will be obtained as consequences of a more precise
fibrewise theorem, which we describe next.

Let $\CC\to N$ be the space of  pairs $(x,\alpha )$
where $x\in M$ and
$\alpha :[0,1]\to N$  is a continuous path
such that $\alpha (0)=f(x)$, fibred over $N$ by projection to
the other endpoint $\alpha (1)$.
We have a homotopy equivalence $\pi :\CC \to M$ given by $\pi (x,\alpha )=x$.

The homotopy Pontryagin-Thom map defined by
the embedding $(e,f) : M\into V\times N$
as described in \cite[Section 6]{coin} (and in \cite{FHT})
is a pointed map
$$
\tilde F :N\times V^+ \to (N\times V^+)\wedge_N\CC_N^{\pi^*\nu}
$$
with compact supports over $N$.
(To be exact,
the space $\CC$ in \cite{coin} is the space of pairs $(x,\beta )$,
where $x\in M$ and $\beta : [0,1] \to V\times N$ is a path starting at
$\beta (0)= (e(x),f(x))$. We omit here the redundant component in the
contractible space $V$.)
We may then form the fibrewise geometric Hopf invariant
$$
h_V(\tilde F) \in{}_{\phantom{Z/}c}^\ZZ\omega^0_N\{ N\times \Sigma S(LV)_+;\,
(N\times LV^+)\wedge_N
(\CC_N^{\pi^*\nu} \wedge_N \CC_N^{\pi^*\nu})\}\, ,
$$
again as a stable $\ZZ$-equivariant map with compact supports over $N$.

Now the fibre product $\CC\times_N\CC$ of pairs $((x,\alpha ),(y,\beta ))$
such that $\alpha (1)=\beta (1)$ may be identified, by splicing
$\alpha$ to the reversed path $\beta$, with
the space $\CD$ of triples $(x,y,\gamma )$ with $x,\, y\in M$ and $\gamma :
[-1,1]\to N$ a continuous path from $\gamma (-1)=f(x)$ to
$\gamma (1)=f(y)$ projecting to $\gamma (0)\in N$.
(Thus $\gamma (t)$ is $\alpha (1+t)$ if $-1\leq t\leq 0$,
$\beta (1-t)$ if $0\leq t\leq 1$.)
It has an action of $\ZZ$ in which the involution interchanges
$x$ and $y$ and reverses the path $\gamma$.
The double point set $\Dp$ is included in $\CD$ as the space of constant paths.

Let $\EE$ be the space of pairs $((x,y),\alpha )$, where
$(x,y)\in \Dp$ and $\alpha :[0,1]\to N$ is a path
such that $\alpha (0)=f(x)=f(y)$.
This corresponds to the subspace of points $(x,y,\gamma )\in \CD$
with $\gamma (-t)=\gamma (t)$.
The fibrewise Pontryagin-Thom construction
on $\Dp\into (LV-\{ 0\})\times N$
gives a fibrewise map
$$
N\times \Sigma S(LV)_+ \to
(N\times LV^+)\wedge_N \EE_N^{\nu\times\nu},
$$
which we compose with the inclusion
$$
\EE_N^{\nu\times\nu}\into \CD_N^{\nu\times\nu},
$$
to get an equivariant fibrewise map
$$
\tilde \phi : N\times \Sigma S(LV)_+ \to
(N\times LV^+)\wedge_N \CD_N^{\nu\times\nu}.
$$

\Thm{\label{hdpt}
{\rm (Homotopy double point theorem).}
The fibrewise geometric Hopf invariant
$$
h_V(\tilde F) \in
{}_{\phantom{Z/}c}^\ZZ\omega^0_N\{ N\times \Sigma S(LV)_+;\,
(N\times LV^+)\wedge_N
\CD_N^{\nu\times\nu}\}
$$
of the homotopy Pontryagin-Thom map $\tilde F$
is equal to the fibrewise stable class of the map $\tilde \phi$
determined by the double points of $f$.
}

The dual version is a class
$$
h'_V(\tilde F) \in
{}_{\phantom{Z/}c}^\ZZ\omega^0_N \{ N\times S^0;\, (N\times S(LV)_+) \wedge_N
\CD_N^{\nu\times\nu}\}\, .
$$
There is also a non-equivariant form.
The stable Hopf invariant $h''_V(\tilde F)$ lies in
$$
{}_c\omega^0_N\{ N\times S^0;\,
(S(LV)\times_{\ZZ} \CD )_N^{\nu\times \nu}\}\, ,
$$
and this group can be identified with
$$
\tilde\omega_0((S(LV)\times_{\ZZ} \CD )^{\nu\times\nu -\tau N})
$$
by fibrewise Poincar\'e-Atiyah duality. (For a general treatment of
fibrewise duality see, for example, \cite[Part II, Section 12]{FHT}.
The duality theorem required
here is stated in \cite{coin} as Proposition 4.1.)

The Pontryagin-Thom construction applied to the
double point manifold $\barDp =\Dp /(\ZZ)$ equipped with the map
$\barDp\to S(LV)\times_{\ZZ} \CD$ given
by the inclusion $\Dp\to\CD$ and the map $\Dp\to S(LV)$ given by $e$
(via $d$) produces a stable homotopy class
$$
\tilde \phi'' :
S^0 \to (S(LV)\times_\ZZ \CD )^{\nu\times\nu -\tau N} .
$$
\Cor{We have
$$
h''_V(\tilde F)=\tilde \phi'' \in
\tilde\omega_0((S(LV)\times_{\ZZ} \CD )^{\nu\times\nu -\tau N})\, .
$$
}
Before turning to the proof of Theorem \ref{hdpt},
we explain how the fibrewise Hopf invariant $h_V(\tilde F)$ determines
the simpler invariant $h_V(F)$.
\Lem{The Hopf invariant $h_V(F)$ is the image of $h_V(\tilde F)$ under
the composition:
$$
\begin{array}{ll}
&{}^\ZZ_{\phantom{Z/}c}\omega^0\{ N\times \Sigma S(LV)_+;\,
(N\times LV^+)\wedge_N \CD^{\nu\times\nu}_N\} \\[1ex]
&\qquad\to
{}^\ZZ\omega^0\{ N^+\wedge \Sigma S(LV)_+;\,
LV^+\wedge \CD^{\nu\times\nu}\}\\[1ex]
&\qquad\to
{}^\ZZ\omega^0\{ N^+\wedge \Sigma S(LV)_+;\,
LV^+\wedge (M^\nu\wedge M^\nu )\}
\end{array}
$$
of the homomorphism defined
by collapsing the basepoint sections over $N$ to a point
and that induced by
the projection $\pi\times\pi :\CD = \CC\times_N\CC\to M\times M$.
}
\begin{proof}
This is easily seen from the explicit construction of the geometric
Hopf invariant.
\end{proof}

In a similar manner, $h'_V(F)$ and $h''_V(F)$ are the images
of the refined invariants $h'_V(\tilde F)$ and $h''_V(\tilde F)$
under homomorphisms defined by collapsing basepoint sections over $N$
and projecting from $\CD$ to $M\times M$.
\Rem{This construction also provides the $\pi$-equivariant Hopf invariant
considered, in greater detail, in \cite{Hopf}.\\
(i) Suppose that $\pi$ is a discrete group and that $q:\widetilde N\to N$
is a principal $\pi$-bundle (for example, a universal covering
space with $\pi$ the fundamental group of $N$).
We let $\widetilde M=f^*\widetilde N$ be the induced bundle over $M$:
thus
$$
\widetilde M=\{ (x,z)\in M\times\widetilde N\st f(x)=q(z)\}\, .
$$
We may define a map $\CD \to (\tilde M\times \tilde M)/\pi$ as follows.
Given $(x,y,\gamma )\in\CD$ (so $x,\, y\in M$, $\gamma : [-1,1]\to N$,
$\gamma (-1)=f(x)$, $\gamma (1)=f(y)$),
lift $\gamma$ to a path $\tilde \gamma$ in $\widetilde N$, determined up
to multiplication by an element of $\pi$. The $\pi$-orbit of
$((x,\tilde\gamma (-1)),(y,\tilde\gamma (1)))\in \widetilde M\times\widetilde M$
is independent of the choice of the lift.
The Hopf invariant $h''_V(\tilde F)$ thus gives us an element of
$$
\omega^0\{ N^+ ;\,
(S(LV)\times_\ZZ ((\widetilde M\times\widetilde M)/\pi ))
^{\nu\times\nu }\}\, ,
$$
where $\nu\times\nu$ is lifted from $M\times M$ to
$(\widetilde M\times\widetilde M)/\pi$ by the obvious projection.\\
(ii) Let $(f,b):N \to X$ be a normal map from an $n$-dimensional manifold
$N$ to an $n$-dimensional geometric Poincar\'e complex $X$.
Let $\widetilde{X}$ be the universal cover of $X$,
let $\widetilde{N}=f^*\widetilde{X}$ be the pullback cover of $N$, and
let $\pi=\pi_1(X)$. The Umkehr $\Zz[\pi]$-module chain map
$$f^!~:~C(\widetilde{X})~\simeq~C(\widetilde{X})^{n-*}
\xymatrix{\ar[r]^-{\widetilde{f}^*}&} C(\widetilde{N})^{n-*}~\simeq~
C(\widetilde{N})$$
is induced by a $\pi$-equivariant geometric Umkehr map
$F:\Sigma^{\infty}\widetilde{X}_+ \to \Sigma^{\infty}\widetilde{N}_+$
$\pi$-equivariant $S$-dual to $T(\widetilde{b}):T(\widetilde{\nu}_N)
\to T(\widetilde{\nu}_X)$. The Wall surgery obstruction $\sigma_*(f,b)
\in L_n(\Zz[\pi])$ was identified in \cite{Ranicki} with the cobordism
class of an $n$-dimensional quadratic Poincar\'e complex $(C(f^!),\psi)$
over $\Zz[\pi]$, with $\psi$ the evaluation of the
`quadratic construction' $\psi_F$ on the fundamental class $[X] \in H_n(X)$.
In \cite{Hopf} we shall express $\psi_F$ in terms of the $\pi$-equivariant
geometric Hopf invariant of $F$.\\
(iii) See \cite{Klein} for an application of the $\pi$-equivariant Hopf invariant
to diagonals in geometric Poincar\'e complexes.
}

In view of the relations between the various Hopf invariants, it
will suffice to prove Theorem \ref{hdpt} in order to establish Theorem \ref{dpt}
and their sundry corollaries.

\medskip

\noindent
{\it Proof of Theorem \ref{hdpt}.}
Writing $\bar \nu$ for the normal bundle of $\barDp$ in $N$, choose
a tubular neighbourhood $D(\bar\nu )\into N$ of $\barDp$.
For $(x,y)\in \Dp$, the fibre of $\bar\nu$ at $f(x)=f(y)$ is
$\nu_x\oplus\nu_y$.
We identify $f^{-1}(\barDp)$ with $\Dp$ by projecting to the first factor.
Then the inverse image of the tubular neighbourhood $D(\bar\nu )$
is a tubular neighbourhood $D(\nu')$ of $\Dp$ in $M$,
where $\nu'_x=\nu_y$. The normal bundle $\nu$ restricted to $D(\nu ')$
is then identified with the restriction of $\nu$ to $\Dp$.

We may assume that $e$ vanishes outside the tubular neighbourhood $D(\nu')$.
The homotopy Pontryagin-Thom map is then a $V$-fold suspension
outside $D(\bar\nu )$. By Proposition \ref{suspension}, the fibrewise
Hopf invariant is canonically null-homotopic outside  $D(\bar\nu )$.
The Hopf invariant $h_V(\tilde F)$ is thus represented by a fibrewise
map which is zero outside the tubular neighbourhood.
In this way we localize $h_V(\tilde F)$ to
an element of
$$
{}^\ZZ\omega^0_{(D(\bar\nu ),S(\bar\nu ))}\{ N\times \Sigma S(LV)_+;
\, (N\times LV^+)\wedge_N \CD^{\nu\times\nu}_N\}
$$
constructed from the immersion data in a neighbourhood of the
double points.

The local data consists simply of the double cover $\Dp\to \barDp$
and the vector bundle $\nu$ over $\Dp$. The bundle $\nu'$ is the
pullback of $\nu$ under the covering involution and $\bar\nu$
over $\barDp$ is the push-forward of $\nu$.
The `local $M$' is the total space of $\nu'$ over $\Dp$, and the `local
$N$' is the total space of $\bar\nu$ over $\barDp$. The immersion
$f$ is given by the projection $\Dp\to\barDp$.
We also need the map $e$ and, without loss of generality, we may assume
that it is determined by a $\ZZ$-map $\Dp\to S(LV)$ (extended radially
on $D(\nu')$ to taper to zero).

In the fibre over $\{ x,y\}\in\barDp$, the manifold $M$ is the disjoint
union
$$
(\{0\}\times\nu_y)\sqcup (\nu_x\times\{ 0\}) \imm\nu_x\oplus\nu_y
$$
immersed in the fibre of $N$ by projection, with the double point at $(0,0)$.
Write $e(x)=v$, $e(y)=-v$. The Hopf invariant is calculated by
the sum formula of Proposition \ref{sum}.
In the notation used there, we take $X= \nu_x^+\times\nu_y^+$,
$Y=\nu_x^+\vee\nu_y^+$,
and the maps $F_+$ and $F_-$ are suspensions of
the compositions of the projection to $\nu_x^+$
or $\nu_y^+$, respectively, and the inclusion of the wedge summand.
The Hopf invariants of the suspensions $F_+$ and $F_-$ vanish
and the Hopf invariant of the sum is determined by the product
$X\to Y\times Y$, so by the projection $\nu_x^+\times\nu_y^+\to
\nu_x^+\wedge \nu_y^+$.

The same computation performed fibrewise over $\barDp$ gives
the localized fibrewise Hopf invariant
as the image of the element in
$$
\omega^0_{(D(\bar \nu ),S(\bar\nu ))}\{ D(\bar\nu )\times S^0;\,
\bar\nu^+_{D(\bar\nu )}\}
$$
given by the inclusion of $D(\bar\nu )$ in $\bar\nu$.
Unravelling the definitions, one sees that this produces $[\phi ]$.
\qed

\smallskip

In the remainder of the paper we shall be concerned with the behaviour
of the geometric Hopf invariant $h_V(\tilde F)$ as one deforms the
map $f: M\to N$.
Suppose, to begin, that we have
smooth homotopies $f_t : M\to N$ and $e_t:M\to V$
such that each $f_t$ is an immersion
and such that $(e_t,f_t) : M\to V\times N$ is an
embedding for each $t\in [0,1]$.
We write $f=f_0$, $f'=f_1$, and $e=e_0$, $e'=e_1$.
Then we have, up to homotopy, an isomorphism
$\nu'=\nu (f_1)\to \nu (f_0)=\nu$
between the normal bundles of the immersions
and a fibre homotopy equivalence $\CDone \to \CD$ over $N$.
In the standard terminology, the homotopy $f_t$ is a {\it regular
homotopy} from $f$ to $f'$ and the homotopy $(e_t,f_t)$ is
an {\it isotopy} from $(e,f)$ to $(e',f')$.
Let $F'$ be the map obtained from the homotopy Pontryagin-Thom
construction on $(e',f')$.
\Prop{{\rm (Regular homotopy/isotopy invariance).}
The fibrewise Hopf invariants $h_V(\tilde F')$ and $h_V(\tilde F)$
correspond under the induced isomorphism:
$$
\begin{array}{ll}
&{}_{\phantom{Z/}c}^\ZZ\omega^0_N\{ N\times \Sigma S(LV)_+;\,
(N\times LV^+)\wedge_N
\CDone_N^{\nu'\times\nu'}\} \\[1ex]
&\to
{}_{\phantom{Z/}c}^\ZZ\omega^0_N\{ N\times \Sigma S(LV)_+;\,
(N\times LV^+)\wedge_N
\CD_N^{\nu \times\nu }\} \, .
\end{array}
$$
}
\begin{proof}
This follows from the homotopy invariance of the geometric Hopf invariant,
which in turn follows easily from the explicit construction described
in the Introduction.
\end{proof}

\Sect{Immersions and embeddings}
The vanishing of the stable Hopf invariant $h_V(\tilde F)$ is a
necessary condition for the existence of a regular homotopy $f_t$
and isotopy $(e_t,f_t)$ such that $f'=f_1$ is an embedding.
In the opposite direction, we record first:
\Prop{Suppose that $m<2(n-1)$ and that $h_V(\tilde F)=0$. Then
$\tilde F$ is
homotopic (through maps with compact support over $N$) to the
$V$-fold suspension
of a map determined by a section $N\to\CC^{\pi^*\nu}_N$.
}
\begin{proof}
This is a consequence of
the fibrewise EHP-sequence (see \cite[Part II, Proposition 14.39]{FHT}),
applicable because $0\leq 3(n-m-1)-(n+0)$.
\end{proof}
To proceed further, we shall assume
that the dimension of $V$ is large: $\dim\, V >2m$.
This guarantees that $M$ can be embedded in $V$, and
we may take the map $e:M\to V$ to be an embedding.
In a metastable range, the fibrewise Hopf invariant of $\tilde F$ is
the precise obstruction to the existence of a regular homotopy
from the immersion $f:M\imm N$ to an embedding.
\Thm{\label{HHQ}
{\rm (Haefliger \cite{Haefliger}, Hatcher-Quinn \cite{Quinn}).}
Suppose that $3m<2(n-1)$ and $\dim V>2m$. Then
$$
h''_V(\tilde F) \in \tilde\omega_0((S(LV)\times_\ZZ \CD)^{\nu\times\nu
-\tau N})
= \tilde\omega_0((E\ZZ \times_\ZZ \CD)^{\nu\times\nu
-\tau N})
$$
vanishes if and only if $f$ is homotopic through immersions to
an embedding of $M$ into $N$.
}

We define $\hDp$ to be the subspace of $\CD$ consisting of the
triples $(x,y,\gamma )$ such that $x\not=y$
and call $\hDp$ the
space of {\it homotopy double points} of the immersion $f$.
The Hopf bundle $E\ZZ\times_\ZZ L$ over $B\ZZ$ will be denoted by $H$.
\Lem{\label{free}
Let $\CCC$ be the complement of
$\hDp$ in $\CD$. Then one has a homotopy
cofibration sequence:
$$
\begin{array}{ll}
&(E\ZZ\times_\ZZ \hDp )^{\nu\times\nu -\tau N}
\to  (E\ZZ\times_\ZZ \CD )^{\nu\times\nu -\tau N} \\ [1ex]
&\qquad\qquad
\to  (E\ZZ\times_\ZZ \CCC )^{H\otimes\tau M+\nu\times\nu -\tau N},
\end{array}
$$
in which the first map is given by the inclusion of $\hDp$ in $\CD$
and the second by the homotopy Pontryagin-Thom construction on the
diagonal submanifold $M$ in $M\times M$.
}

As a space over $M$, $\CCC$ is the pullback by $f : M\to N$ of the free loop
space of $N$, $\map (S^1,N)\to N$, fibred over $N$ by evaluation at
$1\in S^1\, (\subseteq\Cc )$.
\begin{proof}
The vector bundle $L\otimes\tau M$, corresponding to $H\otimes \tau M$
over $B\ZZ$, is the normal bundle of the diagonal
inclusion of $M$ in $M\times M$.
Choose a $\ZZ$-equivariant tubular neighbourhood
$D(L\otimes\tau M)\into M\times M$ of the diagonal
in $M\times M$.
The inclusion of the subspace of $\hDp$ consisting of the triples
$(x,y,\gamma )$ such that $(x,y)\notin B(L\otimes \tau M)$
into $\hDp$ is a homotopy equivalence.

The argument used in the proof of Lemma 6.1 in \cite{coin} then shows
that we have a homotopy cofibre sequence.
\end{proof}
\Cor{\label{isomorphism}
The inclusion induces an isomorphism
$$
\tilde\omega_0((E\ZZ\times_\ZZ \hDp )^{\nu\times\nu -\tau N})
\to
\tilde\omega_0((E\ZZ\times_\ZZ \CD )^{\nu\times\nu -\tau N})
$$
provided that $m<n-1$.
}
\begin{proof}
This follows at once from the long exact sequence of the cofibration.
\end{proof}

The involution on $\hDp$ is free, and hence, writing
$\barhDp$ for the quotient $\hDp /(\ZZ )$, we have an isomorphism
$$
\tilde\omega_0((E\ZZ\times_\ZZ \hDp )^{\nu\times\nu -\tau N})
\to
\tilde\omega_0 (\barhDp^{\nu \times\nu -\tau N})\, .
$$

Until now, we have thought of $\CD$, which arose as the
fibre product $\CC\times_N\CC$, as a space over $N$. It also
fibres over $M\times M$ and the fibrewise space $\CD \to M\times M$
is $\ZZ$-equivariantly locally fibre homotopy trivial.
(Compare \cite[Definition 2.3]{coin}.)
In the same way, $\hDp$ is fibred over the complement $M\times M-\Delta (M)$
and $\barhDp$ is fibred over $(M\times M-\Delta (M))/\ZZ$.

\medskip

\noindent
{\it Proof of Theorem \ref{HHQ}.}
We shall use results and terminology from \cite[Section 7]{coin},
which derive from work of Koschorke \cite{Koschorke} and
Klein and Williams \cite{KleinWilliams1}.

Suppose, first, that $M$ is connected.
Put $\tilde B =M\times M - B(L\otimes \tau M))$; it is a manifold
with a free $\ZZ$-action.
Let $B$ be the manifold $\tilde B/\ZZ$
with boundary $\partial B$ the projective bundle $P(\tau M)$.
Let $E$ be the bundle $(\tilde B \times (N\times N))/\ZZ$ over $B$
and $Z\subseteq E$ the diagonal sub-bundle $B\times N=(\tilde B \times N)/\ZZ$.
The fibrewise normal bundle of the inclusion
$Z\into E$ is $H\otimes \tau N$,
where $H$ is the line bundle over $B$ associated to the double
cover.
The $\ZZ$-equivariant square $f\times f: \tilde B\to N\times N$
defines a section $s$ of $E\to B$, which, if the tubular
neighbourhood is chosen to be sufficiently small, has the property that
$s(x)\notin Z_x$ for $x\in\partial B$.
Together, the fibre bundle $E\to B$, the sub-bundle $Z\to B$ and the section $s$
constitute what is called in \cite{coin} an intersection problem.
In the language used there, $s$ is nowhere null on the boundary
$\partial B$ and the homotopy null-set
 $\hNull (s)$, fibred over $B$,
is easily identified with the restriction of
$\barhDp\to (M\times M-\Delta (M))/\ZZ$.
The inclusion $\hNull (s)\into \barhDp$ is, as we have already
noted in the proof of Lemma \ref{free}, a homotopy equivalence.

If $f$ is an embedding, then the section $s$ is nowhere null.
Now the homotopy Euler index
(\cite[Definition 7.3]{coin})
$$
\hInd (s; \partial B) \in
\tilde\omega_0 (\hNull (s)^{H\otimes\tau N - \tau B)})
$$
is an obstruction to deforming $s$, through a homotopy constant on
$\partial B$, to a section that is nowhere null.
By Proposition 7.4 of \cite{coin}
it is the precise obstruction if $\dim B < 2(\dim N -1)$.
We conclude that if $\hInd (s;\, \partial B)=0$ and $m < n -1$,
then
$$
f\times f : M\times M - B(L\otimes\tau M) \to N\times N
$$
is $\ZZ$-equivariantly homotopic,
by a homotopy that is constant on the boundary $SL\otimes\tau M)$,
to a map into $N\times N - \Delta (N)$.
According to Haefliger \cite[Th\'eor\`eme 2]{Haefliger},
if further $3m<2(n-1)$,
then this is a sufficient condition
for $f:M\imm N$ to be homotopic through immersions to an embedding
of $M$ into $N$.

Thus far, the
argument is taken from \cite[Theorem A.4 and Corollary A.5]{KleinWilliams1}.
We now relate the index $\hInd (s;\partial B)$ to the geometric Hopf invariant.
\Lem{The homotopy Euler index
$$
\hInd (s;\partial B)\in
\tilde\omega_0 (\hNull (s)^{\nu \times\nu -\tau N})
$$
corresponds to the fibrewise Hopf invariant
$$
h_V''(\tilde F)\in\tilde\omega_0((E\ZZ\otimes_{\ZZ} \CD)^{\nu\times\nu
-\tau N})\, .
$$
}
\begin{proof}
The correspondence is made via the homotopy equivalence
$\hNull (s)\into \barhDp$ and the isomorphism from Corollary \ref{isomorphism}.

We can assume that $f$ satisfies the conditions for the homotopy double point
theorem.
Then $h_V''(\tilde F)$ is represented by the double point manifold
$\barDp$.
By Proposition 7.8 of \cite{coin}, the homotopy Euler index
is represented by the null-set $\Null (s)$, which is exactly
$\barDp$.
\end{proof}
Hence the vanishing of $h''_V(\tilde F)$ implies, in the metastable range,
that $f$ is regularly homotopic to an embedding.

This has dealt with the case in which $M$ is connected.
Now suppose that $M$ is a disjoint union $M^{(1)}\sqcup M^{(2)}$
and write $f^{(i)}$ for the restriction of $f$ to $M^{(i)}$.
Then $\CD$ splits equivariantly over
$$
M\times M = (M^{(1)}\times M^{(1)})
\sqcup (M^{(2)}\times M^{(2)}) \sqcup
(M^{(1)}\times M^{(2)}\sqcup M^{(2)}\times M^{(1)})
$$
as a disjoint union
$$
\CDa\sqcup \CDb\sqcup (\ZZ \times \hI )\, ,
$$
where $\hI$ is the space of triples $(x,y,\gamma )$ with
$(x,y)\in M^{(1)}\times M^{(2)}$ and $\gamma (-1)=f^{(1)}(x)$,
$\gamma (1)=f^{(2)}(y)$.
The fibrewise Hopf invariant decomposes, according to
Proposition \ref{sum}, as a sum of three terms.
The first two are the fibrewise Hopf invariants of $f^{(1)}$ and
$f^{(2)}$; the third, $(12)$-component,
is a more elementary product obstruction.

Arguing by induction, we may suppose that both $f^{(1)}$ and $f^{(2)}$
are embeddings, intersecting transversely,
and that $M^{(2)}$ is connected. We consider a new
intersection problem with $B=M^{(2)}$, $E=B\times N$ and
$Z=B\times f(M^{(1)})$. Let $s : B\to E$ be the section
$s(x)=(x,f^{(2)}(x))$.
This time $\hNull (s)$ is $\hI$,
and we may identify the homotopy Euler index $\hInd (s)$ in
the same way with the $(12)$-component of the fibrewise Hopf
invariant, both being represented by the manifold
$f(M^{(1)})\cap f(M^{(2)})$.
By Proposition 7.4 of \cite{coin} again, the vanishing of the
homotopy Euler index implies that $f^{(2)}$ is homotopic to a map
into $N-f(M^{(1)})$, because $\dim M^{(2)} < 2(n -\dim M^{(1)}-1)$.
Finally, we may apply \cite[Theorem 1.1]{Quinn}
to deduce that $f^{(2)}$ is isotopic to an embedding of
$M^{(2)}$ into the complement of $f(M^{(1)})$.
This inductive step is enough to conclude the proof of Theorem \ref{HHQ}.
\qed

\Rem{\label{fun}
Suppose that $M$ (as well as $N$) is connected. Choose basepoints $*\in M$
and $*=f(*)\in N$. Then we can include the loop-space
$\Omega N$ in $\CD$ by mapping
a loop $\gamma : [-1,1]\to N$, with $\gamma (-1)=*=\gamma (1)$,
to $(*,*,\gamma )\in\CD$, and the set of path components of $\CD$ is
in this way identified with the set of double cosets
$$
f_*\pi_1(M)\backslash\pi_1(N)/f_*\pi_1(M) =\pi_0(\CD )
$$
with the $\ZZ$-action given by the group-theoretic inverse. Thus
we may identify
the set of path components of $S(LV)\times_\ZZ \CD$, if $\dim V>1$, with
the orbit space of the involution.
}
\Ex{Suppose that $n=2m$.
Then $\tilde\omega_0((E\ZZ\times_\ZZ \CD )^{\nu\times\nu -\tau N})$
is a direct sum of groups indexed by $\pi_0(\CD )/\ZZ$,
each component being isomorphic to $\Zz$ or $\Zz /2\Zz$.
When $M$ is connected we may label the summands as in Remark \ref{fun}
by equivalence classes of group elements $g\in\pi_1(N)$.
Let $w_M : \pi_1(M)\to \{ \pm 1\}$ and $w_N : \pi_1(N)\to \{\pm 1\}$
be the orientation maps (corresponding to $w_1M$ and $w_1N$).
The $g$-summand is isomorphic to $\Zz$ if and only if
for all $a,\, b\in\pi_1(M)$:
(o) if $f_*(a)g=gf_*(b)$, then $w_M(ab)=w_N(f_*(a))$ ($=w_N(f_*(b))$),
and
(i) if $f_*(a)g=g^{-1}f_*(b)$, then $(-1)^mw_M(ab)=w_N(f_*(a)g)$.
In particular, if $M$ is orientable,
$N$ is oriented and $f_*\pi_1(M)$ is trivial,
then the obstruction group is
$\Zz [\pi_1(N)]/\langle g-(-1)^mg^{-1}\st g\in\pi_1(N)\rangle$
and the Hopf invariant $h''_V(\tilde F)$
is Wall's invariant $\mu (f)$ in the form mentioned in the Introduction.
}
\Sect{Homotopic immersions}
We next investigate immersions homotopic, as maps, to the given immersion
$f$.
Consider smooth homotopies $f_t : M\to N$ and $e_t : M\to V$
such that $f_0$ and $f_1$ are immersions and each map
$(e_t,f_t):M\to V\times N$, for $0\leq t\leq 1$, is an embedding with
normal bundle $\mu_t$.
We again write $f=f_0$, $f'=f_1$, and $\nu =\nu (f)$, $\nu'=\nu (f')$.
The homotopies determine, up to a homotopy,
a bundle isomorphism $a:V\oplus \nu'=\mu_1\to\mu_0 =V\oplus \nu$
and a $\ZZ$-equivariant fibre homotopy equivalence $\CDone \to \CD$ over
$M\times M$.

There is an associated class
$$
\theta (e_t,f_t)\in\tilde\omega_0((P(V)\times M)^{H\otimes\nu -\tau M})
$$
which vanishes if each $f_t$ is an immersion.
It is constructed as follows.
Let $v_0 : V\to V\oplus \nu$ be the inclusion and
let $v_1: V\to V\oplus \nu'=\mu_1 \to\mu_0=V\oplus\nu$
be the composition of the inclusion with the isomorphism $a$.
Now we have a stable cohomotopy difference class
$$
\delta (v_0,v_1)\in
\tilde\omega^{-1}((M\times P(V))^{-H\otimes (V\oplus\nu )})\, ,
$$
which is the metastable obstruction to deforming $v_0$ to $v_1$
through vector bundle monomorphisms;
see Section 6 for details.
We define $\theta (e_t,f_t)$ to be the dual class in stable homotopy.
(Recall that $\tau P(V)\oplus\Rr =H\otimes V$.)

Using the inclusion $M\into\CD$ : $x\mapsto (x,x,\gamma )$,
where $\gamma$ is the constant loop at $f(x)$, we get a map
$$
i:
\tilde\omega_0((P(V)\times M)^{H\otimes \nu -\tau M})
\to
\tilde\omega_0((S(LV)\times_{\ZZ} \CD )^{\nu\times\nu-\tau N})\, .
$$
\Rem{The map $M\to\CD$ picks out
(under the assumption that $M$ is connected)
a component $\CD_0$ of $\CD$,
which is preserved by $\ZZ$.
Hence $i$ maps into the summand
$$
\tilde\omega_0((S(LV)\times_\ZZ \CD_0)^{\nu\times\nu -\tau N})\, .
$$
}

The link between the difference class $\theta$ and the Hopf invariant
is forged by a generalization of the classical description of
the Hopf invariant on the image of the $J$-homomorphism.
\Prop{\label{ImJ}
Let $a : V\oplus \nu' \to V\oplus \nu$ be a vector bundle isomorphism
over $M$. The fibrewise one-point compactification of $a$ is a map
of sphere bundles
$$
A: (M\times V^+)\wedge_M (\nu')^+_M \to (M\times V)^+\wedge \nu^+_M
$$
over $M$. Then the fibrewise Hopf invariant
$$
h_V(A) \in {}^\ZZ\omega^0_M\{ (M\times \Sigma S(LV)_+)\wedge_M (\nu')^+_M
;\, (M\times LV^+)\wedge_M (\nu^+_M \wedge_M \nu^+_M) \}
$$
of $A$
coincides up to sign,
under the identifications described below, with the difference class
$$
\delta (v_0,v_1) \in
\tilde\omega^{-1}((M\times P(V))^{-H\otimes(V\oplus \nu )}))
$$
of the monomorphisms $v_0,\, v_1 : V\to V\oplus\nu$ given, respectively,
by the inclusion of the first factor and the composition of the inclusion
$V\to V\oplus\nu'$ with $a$.
}
\begin{proof}
The fibrewise smash product $\nu^+_M\wedge_M\nu^+_M
=(\nu\oplus\nu )_M^+$ with the action of $\ZZ$ which
interchanges the factors is equivariantly homeomorphic
(by the construction $\kappa_V$ in Section 1) to
$(\nu\oplus L\nu )^+_M = \nu^+_M\wedge (L\nu )^+_M$. The
isomorphism $a: V\oplus \nu'\to V\oplus\nu$ gives a stable fibre
homotopy equivalence $(\nu')^+_M\to\nu^+_M$. Taken together,
these equivalences allow us to think of $h_V(A)$ as an element of
the group
$$
{}^\ZZ\omega^0_M\{ M\times \Sigma S(LV)_+;\,
(LV\oplus L\nu )^+_M\}\, ,
$$
which is then identified with
$$
{}^\ZZ\tilde\omega^{-1}((M\times S(LV))^{-(LV\oplus L\nu)})
=
\tilde\omega^{-1}((M\times P(V))^{-H\otimes(V\oplus \nu )}) \, .
$$
Both $h_V(A)$ and $\delta (v_0,v_1)$ are defined by difference constructions.
The proof that they coincide follows from
a direct comparison of the definitions.
\end{proof}
\Thm{\label{Var}
{\rm (Homotopy/isotopy variation).}
The Hopf invariant
$$
h''_V(\tilde F')\in \tilde\omega_0((S(LV)\times_\ZZ \CDone)
^{\nu'\times\nu' -\tau N})
$$
corresponds to
$$
h''_V(\tilde F) + i\theta (e_t,f_t) \in \tilde\omega_0((S(LV)\times_\ZZ \CD)
^{\nu \times\nu  -\tau N})\, .
$$
}
\begin{proof}
Recollect that $\tilde F$ is a fibrewise map with compact supports
over $N$:
$$
N\times V^+ \to (N\times V^+)\wedge_N \CC_N^{\pi^*\nu}.
$$
Up to homotopy, the map $\tilde F'$ associated with $(e',f')$ is
the composition of $\tilde F$ with the map
$$
(N\times V^+)\wedge_N \CC_N^{\pi^*\nu} =\CC_N^{\pi^*(V\oplus\nu )} \to
\CC_N^{\pi^*(V\oplus \nu')}=(N\times V^+)\wedge_N \CC_N^{\pi^*\nu'}
$$
obtained by lifting the bundle isomorphism $a^{-1} : V\oplus \nu
\to V\oplus \nu'$ over $M$ via $\pi : \CC\to M$.
In the notation of Proposition \ref{ImJ} we have $\tilde F\simeq
\pi^*A\comp\tilde F'$.

The fibrewise version of the composition formula
(Proposition \ref{composition}) expresses the
difference $\alpha =h_V(\tilde F)-h_V(\tilde F')$
in terms of $A$.
Using Proposition \ref{ImJ} and the explicit form of the geometric Hopf
invariant one sees that $\alpha$ lies in the image of the diagonal
map
$$
\begin{array}{ll}
\Delta_* :
&{}^\ZZ_{\phantom{Z/}c}\omega_N^0\{ N\times \Sigma S(LV)_+;\,
(N\times LV^+)\wedge_N \CC^{\pi^*(\nu \oplus L\nu )}_N\} \\
&\qquad
\to {}^\ZZ_{\phantom{Z/}c}\omega_N^0\{ N\times \Sigma S(LV)_+;\, (N\times LV^+)\wedge_N
(\CC^{\pi^*\nu}_N\wedge_N\CC^{\pi^*\nu )}_N)\} \, .
\end{array}
$$
This may be rewritten in dual form as:
$$
\tilde\omega_0((P(V)\times \CC)^{\pi^*(\nu \oplus H\otimes\nu -\tau N)})
\to \tilde\omega_0((S(LV)\times_\ZZ \CD) ^{\nu \times\nu  -\tau N})\, .
$$
But $\pi :\CC \to M$ is a homotopy equivalence.
Hence $\Delta_*$ is just another manifestation of the map $i$ in
the statement of the theorem.

Now the tubular neighbourhood of $M$ in $V\times N$ gives a proper
map $B(\nu ) \to N$. To calculate $\alpha$, which is concentrated on
the diagonal, we can thus lift from $N$ to $B(\nu )$.
Here we have a fibrewise problem over $M$, and the identification
of $\alpha$ is achieved by Proposition \ref{ImJ}.
\end{proof}
\Rem{In the stable range $\dim V >2m$, where the maps $e_t$ are
redundant, we can use the methods of
the previous section to give an alternative proof of Theorem \ref{Var}.
}
\Thm{{\rm (Hirsch \cite{Hirsch}).}
Suppose that $3m<2n-1$ and $\dim V >2m$.
Then two immersions $f$ and $f'$ are homotopic through immersions
if and only if the associated difference class
$\theta$ in $\tilde\omega_0((P(V)\times M)^{H\otimes\nu -\tau M})$
is zero.
}
\begin{proof}
The derivative of the immersion $f$ gives a vector bundle monomorphism
$df : \tau M \to f^*\tau N$ over $M$.
According to Hirsch \cite{Hirsch}, for $m<n$
immersions $f': M\to N$
together with a homotopy $f_t$ from $f=f_0$ to $f'=f_1$
are classified by homotopy
classes of vector bundle monomorphisms $\tau M \to f^*\tau N$.
In the metastable range $m+1 < 2(n-m)$, that is, $3m<2n-1$,
immersions with a homotopy to $f$ are thus classified by
$$
\tilde\omega^{-1}(P(\tau M)^{-H\otimes f^*\tau N})
=
\tilde\omega_0(P(\tau M)^{H\otimes ( f^*\tau N-\tau M)-\tau M})
=
\tilde\omega_0(P(\tau M)^{H\otimes \nu-\tau M})\, .
$$
Assuming that $\dim V>2m$ we may fix an embedding $e : M\into V$.
The derivative of $e$ includes
$\tau M$ in the trivial bundle $M\times V$ and gives an isomorphism
$$
\tilde\omega_0(P(\tau M)^{H\otimes \nu -\tau M}) \to
\tilde\omega_0((M\times P(V))^{H\otimes \nu -\tau M})\, .
$$
\vskip-1.4\baselineskip
\end{proof}

\Sect{Immersions close to an embedding}
Consider the special case of a closed manifold $M$ of
(constant) dimension $m$
and a real vector bundle $\nu$ of dimension $n-m$ over $M$.
Working in the metastable range $3m<2n-1$,
we take $N$ to be the total space of $\nu$ and $f: M\to N$ to be the
embedding given by the zero section of the vector bundle.
As $e:M\to V$ we may take the constant map $0$.

Let $v_0 : V\into V\oplus\nu$ be the inclusion of the first factor.
Suppose that $v_1: V\into V\oplus \nu$ is another inclusion,
which we may assume to be isometric.
Thus $v_0$ and $v_1$ give sections of the bundle $\O(V,V\oplus \nu )$
whose fibre at $x\in M$
is the Stiefel manifold of orthogonal linear maps
$v:V\to V\oplus\nu_x$.
Let $X_1(V,\nu )$ be the sub-bundle with fibre consisting of those
linear maps $v$ such that $v+i_x$ has kernel of dimension $1$,
where $i_x:V\to V\oplus\nu_x$ is the inclusion of the
first factor.
In the range of dimensions that we are considering,
for a generic smooth section $v_1$ of $\O(V,V\oplus\nu )\to M$
the kernel of $(v_1)_x+i_x$ is nowhere of dimension $>1$,
$v_1$ is transverse to the sub-bundle $X_1(V,\nu )$ and
$v_1^{-1}(X_1(V,\nu ))$ is a submanifold
$\barD$ of $M$ of dimension $2m-n$, equipped with a map $\barD\to P(V)$
classifying the $1$-dimensional kernel.
The normal bundle of $\barD$ in $M$ is identified with
$H\otimes\nu$ and $\barD$ represents the element
$\delta\in\tilde\omega_0((M\times P(V))^{H\otimes\nu -\tau M})$
dual to $\delta (v_0,v_1)$.
More details are provided in an appendix (Section 6).

We want to construct an immersion $f'$ close to the zero section $f$
with double point set $\barD$. More precisely, we shall construct
homotopies $f_t$ and $e_t$, with $f_0=f$, $f_1=f'$ and each $(e_t,f_t)$
an embedding, such that $\barDpone =\barD$ and $\theta (e_t,f_t)=\delta$.

Choose an open tubular neighbourhood $\Omega =H\otimes\nu\into M$ of $\barD$.
Since $\dim \nu =n-m >2m-n=\dim \barD$, we can split
off a trivial line from the restriction of $\nu$ to $\barD$ as:
$\nu |\barD =\Rr\oplus\zeta$.

Whitney gave in \cite{Whitney}, for a Euclidean space $U$,
an explicit `punctured figure-of-eight' immersion:
$$
w: \Rr\oplus U\to (\Rr\oplus U)\times (\Rr\oplus U)
$$
with double points at $(\pm 1,0)$.
In slightly modified form it may be written as
$$
w(s,y) = ((1-\lambda (s,y))s,y,-\lambda (s,y),\lambda (s,y) sy)\, ,
$$
where $\lambda (s,y)=\psi (s^2+\| y\|^2)$
and $\psi : [0,\infty )\to \Rr$ is
a smooth, non-negative, monotonic decreasing function,
with $\psi (1)=1$, $\psi'(1)=-1/2$ and $\psi (r)=0$ for $r\geq 2$.
(In \cite{Whitney}, $\psi (r)=2/(1+r)$.)
The derivative at the double point $(\pm 1,0)$ is
$$
{\textstyle
\dfrac{\partial w}{\partial s} = (1,0,\pm 1,0),\qquad
\dfrac{\partial w}{\partial y} = (0,1,0,\pm 1).
}
$$
Writing $w(s,y)$ in the form
$(s,y,0,0) +\lambda (s,y)(-s,0,-1,sy)$, we see that
$w(s,y) =(s,y,0,0)$ for $\|(s,y)\|$ large.
The immersion $w$ has $\ZZ$-symmetry as an equivariant map
$$
w: L\oplus LU\to (L\oplus LU)\times (\Rr\oplus U)\, .
$$
The two double points are distinguished by the $\ZZ$-map
$$
c : L\oplus LU \to L
$$
given by $c(s,y)= \lambda (s,y)s$.

Now Whitney's construction, applied on the fibres of $\zeta$,
gives an immersion
$$
\begin{array}{ll}
f': &\Omega = H\otimes\nu =H\oplus (H\otimes\zeta ) \\
&\qquad \to
(H\oplus (H\otimes\zeta ))\times (\Rr\oplus\zeta) =
(H\otimes\nu )\times\nu  = \nu\, |\, \Omega \subseteq N
\end{array}
$$
of the open subset $\Omega$ of $M$ into the total space of the restriction
of $\nu$ to $\Omega$.
We extend $f'$ to the whole of $M$ to coincide with the zero section
$f$ outside a compact subspace of $\Omega$.
Its double point set $\Dpone$ is the double-cover $S(H|\barD )$ of $\barD$ in
$\Omega$.
The map $c$ composed with the inclusion of $H$ into the trivial bundle
$\Omega\times V$ and the projection to $V$ gives a map
$$
e' : \Omega = H\oplus (H\otimes\zeta ) \to V\, ,
$$
which is zero outside a compact subset of $\Omega$ and can be extended
by $0$ to a map $e': M\to V$ which distinguishes the double point
pairs.
The required homotopies $e_t$ and $f_t$, $t\in [0,1]$, joining
$e$ to $e'$ and $f$ to $f'$ are
defined by replacing $\lambda$ in the definition of $e'$ and $f'$
by $t\lambda$.
One checks that $(e_t,f_t)$ is an embedding for all $t$, but that
$f_t$ fails to be an immersion when $t=1/\psi (0)$.
\Thm{Suppose that $3m<2n-1$.
Then Whitney's construction described above produces, for any given element
$\delta\in\tilde\omega_0((M\times P(V))^{H\otimes\nu -\tau M})$, a
homotopy $(e_t,f_t)$ with $\theta (e_t,f_t)=\delta$.
}
\begin{proof}
By construction, the double point manifold $\Dpone =S(H|\barD)$
represents $\delta$.
The assertion thus follows from
the Double Point Theorem \ref{hdpt} for $(e',f')$
in conjunction with the Homotopy Variation Theorem \ref{Var}.
\end{proof}
\Rem{The same construction may be used to modify a general immersion
$f:M\imm N$, and map $e:M\to V$, in the complement of $\Dp$.
We can insert $\barD$ in the complement, because $2(2m-n)<m$.
}
\Ex{Whitney's construction gives the classical immersions of the
sphere $S^m\imm S^{2m}$ close to the equatorial inclusion
and $S^m\imm S^m\times S^m$ close to the diagonal.
}

The early work of Smale \cite{Smale,Smale2} has been followed
by a vast literature on the homotopy-theoretic properties
of immersions, including \cite{Dax}, \cite{HaefligerHirsch}, \cite{Quinn},
\cite{Hirsch}, \cite{KleinWilliams2}, \cite{Ping} and \cite{Salomonsen}.
\Sect{Appendix: Monomorphisms of vector bundles}
Let $\xi$ and $\eta$ be smooth real vector bundles, of dimension
$n$ and $r$ respectively, over a closed
$m$-manifold $M$. We shall describe the differential-topological
classification of homotopy classes of
vector bundle monomorphisms $\eta\to\xi$ in the
metastable range $m +1 < 2(n - r)$.

Suppose that $v_0,\, v_1: \eta \into\xi$ are two vector bundle
monomorphisms.
Doing homotopy theory,
we may assume that $\xi$ and $\eta$ have positive-definite inner
products and that the monomorphisms are isometric embeddings.
Then $v_0$ and $v_1$ are sections of
the bundle $\O (\eta ,\xi )$ whose fibre at $x\in M$ is
the Stiefel manifold of orthogonal linear maps $v: \eta_x\to\xi_x$.
Topological obstruction theory gives a difference class
$$
\delta (v_0 ,v_1) \in \tilde\omega^{-1}(P(\eta )^{-H\otimes \xi})\, ,
$$
where $P(\eta )$ is the projective bundle of $\eta$ and $H$ is
the Hopf line bundle.
This arises as follows. A section of $\O (\eta ,\xi )$
determines a nowhere zero
section of $H\otimes \xi$ over $P(\eta )$:
over $\ell\in P(\eta_x)$ (where $\ell\subseteq \eta_x$ is a line)
a monomorphism $v:\eta_x\to \xi_x$ gives an embedding of
$\ell$ in $\xi_x$ and so a non-zero
vector in $\ell^*\otimes \xi_x$, which is the fibre of
$H\otimes \xi$.
Then $\delta (v_0, v_1)$ is defined as the difference class
$\delta (s_0, s_1)$ of the two nowhere zero sections $s_0$
and $s_1$ of the vector bundle $\xi$ over $P(\eta )$ constructed in this
way from $v_0$ and $v_1$.
We may assume that $s_0$ and $s_1$ are sections of the sphere bundle
$S(H\otimes \xi )$.
Write $s_t = (1-t)s_0 +ts_1$ for $0\leq t\leq 1$. Then $\delta (s_0,s_1)$
is represented explicitly by the map, over $P(\eta )$,
$$
\bar s: ([0,1],\partial [0,1]) \times P(\eta )
\to (D(H\otimes \xi ),S(H\otimes \xi ))
$$
given by the homotopy $s_t$.
In the metastable range, the vector bundle monomorphisms $v_0$ and $v_1$
are homotopic if and only if $\delta (v_0,v_1)=0$.
(See, for example, \cite{Z2, FHT, coin}.)

Thus far, the theory is topological.
We now use Poincar\'e duality for the manifold $P(\eta )$ to identify
$\tilde\omega^{-1}(P(\eta )^{-H\otimes \xi} )$
with $\tilde\omega_0(P(\eta )^{H\otimes (\xi -\eta ) -\tau M})$.
(Up to homotopy, the stable tangent bundle is given by
an isomorphism $\Rr\oplus\tau P(\eta )\iso (H\otimes\eta )\oplus\tau M$.)
Assuming that the monomorphisms $v_0$ and $v_1$ are smooth we shall
represent the dual obstruction class by
a submanifold $Z$ of $M$ together with a map $Z\to P(\eta )$ and appropriate
normal bundle information.
The monomorphism $v_0$ will play a special r\^ole in the description;
to emphasize this,
we write $i=v_0$ for the preferred embedding $i:\eta\into\xi$
and write $\nu$ for the orthogonal complement of $i(\eta )$ in $\xi$.
Let $X_k(\eta ,\nu )$, for $k\geq 1$,
be the sub-bundle of $\O (\eta ,\xi )=\O (\eta ,\eta\oplus\nu )$
with fibre consisting of those linear maps
$v$ such that $i_x+v$ has kernel of dimension $k$.
By Lemma \ref{Cayley} below, we may assume, if $m+1<2(n-r)$, that
$v_1$ never meets $X_k(\eta ,\nu )$ for
$k>1$ and is transverse to the sub-bundle
$X_1(\eta ,\nu )$. The inverse image $v_1^{-1}(X_1(\eta ,\nu ))$
is, therefore, a submanifold
$Z$ of $M$ of dimension $m+r-n$, equipped with a section $Z\to P(\eta )$
classifying the $1$-dimensional kernel. The normal bundle of
$Z$ in $M$ is identified, by Lemma \ref{Cayley}, with $H\otimes\nu$.
\Prop{\label{obstruction}
The submanifold $Z$ described above, with the line bundle $H$
classified by the section of $P(\eta )$ over $Z$ and the isomorphism
between the normal bundle and $H\otimes\nu$, represents the dual
of $\delta (v_0,v_1)$ in  $\tilde\omega_0(P(\eta )^{H\otimes \nu -\tau M})$.
}
\begin{proof}
Consider the sections $s_0$ and $s_1$ of $S(H\otimes\xi )$ over $P(\eta )$
associated with $v_0$ and $v_1$ as in the definition of
$\delta (v_0,v_1)$. The section $\bar s$ of $D(H\otimes\xi )$ over
$[0,1]\times P(\eta )$ given by the homotopy $s_t =(1-t)s_0+ts_1$
is transverse to the zero section and its zero-set is precisely
$\{ \frac{1}{2}\}\times Z$.
The normal bundle is $\Rr \oplus \tau_MP(\eta )\oplus (H\otimes\nu )$,
where $\tau_MP(\eta )$ is the tangent bundle along the fibres of
$P(\eta )\to M$, and this is identified with $(H\otimes\eta )\oplus
(H\otimes\nu )=H\otimes\xi$.
Hence, $Z$ with the normal bundle data represents the stable homotopy class
dual to $\delta (s_0,s_1)=\delta (v_0,v_1)$.
(This is the classical representation of the dual Euler
class of a vector bundle by the zero-set of a generic smooth section.)
\end{proof}
\Rem{A more symmetric treatment may be given by looking
at sections of the fibre product $\O (\eta,\xi )\times_M\O (\eta ,\xi )$
and the sub-bundles with fibre consisting of the pairs
$(u,v)$ such that $\dim\ker (u+v)=k$.
}

The properties of $v_1$ required in the proof of Proposition \ref{obstruction}
follow from  the next lemma, in which
the Lie algebra of the orthogonal group $\O (V)$ of a Euclidean
vector space $V$ is written as ${\mathfrak o}(V)$.
\Lem{\label{Cayley}
Let $V$ and $W$ be finite-dimensional orthogonal vector spaces.
For $0\leq k\leq\dim V$, let $X_k(V,W)$  be the subspace
of the Stiefel manifold $\O (V,V\oplus W)$ consisting of
the maps $v$ such that $i+v$, where $i$ is the inclusion of the first
summand $V\into V\oplus W$, has kernel of dimension $k$.
Then $X_k(V,W)$ is a submanifold diffeomorphic to the total space of
the vector bundle
${\mathfrak o}(\zeta^\perp )\oplus\Hom (\zeta^\perp, W)$
over the Grassmann manifold
$G_k(V)$ of $k$-planes in $V$, where $\zeta^\perp$ is the
orthogonal complement in $V$ of the
canonical $k$-dimensional vector bundle $\zeta$ over $G_k(V)$.
Its normal bundle is naturally identified with
${\mathfrak o}(\zeta )\oplus \Hom (\zeta ,W)$.
}
\begin{proof}
This can be established by using
the (generalized) Cayley transform,
which is written down explicitly in \cite[Part II, Lemma 13.13]{FHT}.
The restriction of the normal bundle of the embedded submanifold
to the subspace $G_k(V)$ is naturally identified with
${\mathfrak o}(\zeta )\oplus \Hom (\zeta ,W)$.
The normal bundle itself is naturally identified with
the pullback by parallel translation.
\end{proof}

In particular, the submanifold $X_k(\eta ,\nu )$ considered above
has codimension $(n-r)k+k(k-1)/2$ in $\O (\eta
,\eta\oplus\nu )$. So, if $k\geq 2$, the codimension is at least
$2(n-r)+1$. The condition $m<2(n-r)+1$ ensures that a generic
section of $\O (\eta ,\eta\oplus\nu )$ is transverse to $X_1(\eta
,\nu )$ and disjoint from $X_k(\eta ,\nu )$ for $k>1$. \Rem{In
\cite{LNM} Koschorke gave an intermediate representation of the
dual of $\delta (v_0,v_1)$ by a submanifold of $(0,1)\times M$.
Consider the section $\bar v$ of $\Hom (\eta ,\xi )$ over
$[0,1]\times M$ given by the homotopy $v_t=(1-t)v_0 +tv_1$ Let
$Y_k(\eta ,\xi )$ be the sub-bundle of $\Hom (\eta, \xi )$ with
fibre consisting of the linear maps with kernel of dimension $k$;
it has codimension $(n-r)k+k^2$, which is $\geq 2(n-r)+4 >m+1$.
Suppose that $\bar v$  can be deformed, by a homotopy through maps
coinciding with $v_0$ and $v_1$ at the endpoints, to a smooth
section $\bar v'$, that never meets $Y_k(\eta ,\xi )$ for $k>1$
and meets $Y_1(\eta ,\xi )$ transversely. This is always possible
if $m<2(n-r)+3$. The inverse image of $Y_1(\eta ,\xi )$ is then a
submanifold $Z$ of $(0,1)\times M$ of dimension $m+r-n$ equipped
with a map $Z\to P(\eta )$ given by the $1$-dimensional kernel.
This data, too, represents the dual of $\delta (v_0,v_1)$. }

\end{document}